\documentclass[10pt,twoside]{article}
\usepackage{amsmath}
\usepackage{amsfonts}

\def  \E {{\mathbb E}}
\def \P {{\mathbb P}}
\def \W {{\cal W}}
\def \D {\cal D} 
\def \S {{\bf S}}

\def \W{{\cal W}}
\def \R{{\bf R}}
\def \B{{\cal B}}
\def \D{{\cal D}}

\def \S{{\bf S}}

\def \e{{\sf e}}

\def \vep {{\varepsilon}}

\def \dd {\xrightarrow[n]{(d)}}

\def \dd1{\xrightarrow[n]{}}

\def \CF {{\cal F}}

\def \L {{ {[\hskip-0.14cm [ \, }}}
\def \R {{{\, ]\hskip-0.14cm]}}}

\begin{document}

\title{ Uniform bounds for exponential moment of maximum of  a Dyck path}

\author{O. Khorunzhiy\footnote{Universit\'e de Versailles, FRANCE}  \  and 
J.-F. Marckert\footnote{CNRS LaBRI, Universit\'e de Bordeaux, FRANCE} }

\maketitle

\begin{abstract}

Let $D_{2n}$ be a Dyck path chosen uniformly in the set of Dyck paths with $2n$ steps.  
The aim of this note is to show that for any $\lambda>0$ the sequence 
${\mathbb E} (\exp(\lambda (2n)^{-1/2}\max D_{2n}))$ converges, and therefore is bounded uniformly in $n$. The uniform bound justifies an assumption used 
in  literature to prove certain \mbox{estimates} of high moments of large random matrices.  
\end{abstract}

\section{Introduction}
Let ${\mathbb N}=\{0,1,2,3,\dots\}$ be the set of non-negative integers.
For any $n\in {\mathbb N}$, we denote by $\W_n$ the set of Bernoulli chains with $n$ steps~:
\[\W_n=\{ {\S}=(S_i)_{0\leq i \leq n}: \ S_0=0, S_{i+1}=S_i \pm 1 \textrm{ for any }i\in
\L{0,n-1}\R \},\]
where $\L{a,b}\R=[a,b]\cap {\mathbb N}$.
The set of Dyck paths $\D_n$ (sometimes called simple or Bernoulli excursions) is defined by
\[\D_n=\{\S~:\S\in\W_n, S_n=0, S_i\geq 0 \textrm{ for any }i \in\L{0,n}\R\}.\]
It is clear that $\D_n$ is empty for odd $n$, and one has
\begin{equation}\label{Cat}
\#\D_{2n}=\frac{1}{n+1}{\binom{2n}{n}}
\end{equation}
 the $n$th Catalan number.
Let $\P^{(w)}_{n}$ and  $\P^{(d)}_{2n}$ be the uniform distributions on $\W_n$ and  $\D_{2n}$. The expectations with respect to these measures will be denoted by $\E^{(w)}_{n}$ and $\E^{(d)}_{2n}$.  The aim of this note is to prove the following statement.
\vskip 0.5cm 
\noindent {\bf Theorem 1} 
\label{princ}
{\it For any $\lambda>0$, we have 
\begin{equation}\label{cvprinc}
\E^{(d)}_{2n}\left(\exp\left(\lambda \frac{\max \S}{\sqrt{2n}}\right)\right)\dd1 
\E\left(\exp\left(\lambda \max_{t \in [0,1]} \e(t)\right)\right)
\end{equation}
where $(\e(t),t\in[0,1]))$ is the normalized Brownian excursion. In particular, for any $\lambda>0$
\begin{equation}\label{pfou}
\sup_{n}\E^{(d)}_{2n}(\exp(\lambda n^{-1/2}\max \S))<+\infty.
\end{equation}
}
\vskip 0.5cm 
We may notice that the right hand side of (2) is finite for every $\lambda$ by using the computation of Chung \cite{CHU}:
\[\P\left(\max _{t \in [0,1]}\e(t) \leq x\right)=1+2\sum_{j\geq 1}(1-4j^2x^2)e^{-2j^2x^2},~~~\textrm{ for }x>0,\]
and then $\max_{t \in [0,1]} \e(t)$ possesses all exponential moments. \par

\subsection{Relations with the spectral theory of large random matrices}

Dyck paths play a central role in combinatorics and arise in many situations (see for instance 66 examples of the appearance of the Catalan numbers in combinatorics in Stanley \cite[ex. 6.19 p.219]{STAN}). 
In the present case, the motivation for establishing Theorem \ref{princ} comes from the study of the moments of  random real symmetric (or hermitian) matrices known as the Wigner ensemble \cite{SS,W}. In these settings, the Catalan numbers (1) 
represent the moments of the eigenvalue distribution of random matrices of Wigner ensemble $A_N$ in the limit of their infinite dimension $N\to\infty$ \cite{W}.
\vskip 0.2cm 
 Recent studies \cite{SS,S}
of high moments of large Wigner random matrices have used the exponential moments 
 of the maxima of the Dyck paths. More precisely, it was shown that \cite{KV,S}
 \begin{equation}\label{boundQ}
 \E\left( \hbox{Tr } A_N^{2n }\right) \le C_1{e^{C_2t^{1/3}} \over t^{1/6}} Q_n (C_3 t^{1/6}) (1+o(1)), \quad n =
  \lfloor tN^{2/3}\rfloor, \ \ N\to\infty, 
 \end{equation}
 where Tr denotes the trace of the square matrix, $C_1$, $C_2$ and $C_3$ are certain constants, and 
\[Q_n(\lambda)=\E^{(d)}_{2n}\left(\exp\left(\lambda \frac{\max \S}{\sqrt{2n}}\right)\right).\]
It was assumed in \cite{KV,SS,S} that $\limsup_n Q_n(\lambda)$ is bounded. 
Theorem 1 above shows that this assertion is  true.

\section{Proof of the Theorem}
Before proving the Theorem, we first discuss the appearance of $\max_{t \in[0,1]} \e(t)$ and the non-triviality of the result. 
Let $C[0,1]$ be the set of continuous functions defined on $[0,1]$ with real values. For any $\S\in \W_n$, denote by 
$u_n=u_n^\S$ the function in $C[0,1]$ obtained from $\S$ by interpolation and rescaling: 
\begin{equation}\label{uuuu}
u_n(t)=\frac{1}{\sqrt{n}}\big(S({\lfloor{nt}\rfloor})+\{nt\}(S({\lceil{nt}\rceil})-S({\lfloor{nt}\rfloor}))\, \big)\textrm{ for any }t\in[0,1].
\end{equation}
It is known that under $\P^{(d)}_{2n}$,  $u_{2n} \xrightarrow[n]{(d)}  \e$ in $C[0,1]$  endowed with the topology of uniform convergence where 
$ \xrightarrow[n]{(d)}$ 
means the convergence in distribution (see e.g. Kaigh \cite{KAI} where this results is shown for general increment distributions). 
By continuity of the map $f \mapsto e^{\lambda \max f}$ from $C[0,1]$ into ${\mathbb R}$,
under $\P^{(d)}_{2n}$, 
\begin{equation}\label{voui}
\exp\left(\frac{\lambda\max \S}{\sqrt{2n}}\right)=\exp(\lambda\max u_{2n})
\xrightarrow[n]{(d)} \exp\left(\lambda\max_{t\in[0,1]} \e(t)\right).
\end{equation}
Then   the uniform integrability argument is sufficient: given $\lambda>0$, in order to prove that (6) implies (2), 
it suffices to show that  
\begin{equation}
\sup_n \E^{(d)}_{2n}\left(\exp\left((\lambda+\vep ) \frac{\max \S}{\sqrt{2n}}\right)\right)<+\infty,
\end{equation}
for some $\vep>0$ (see Billingsley \cite[Section 16]{BIL}).
Hence to prove the Theorem, using 
(6), 
only the second assertion (3) (weaker in appearance) 
needs to be proved. This is what we will do.\par
\vskip 0.3cm 
\noindent {\bf Remark.} 
{\it Smith and Diaconis \cite{SD} proved that 
$\P^{(d)}_{2n}(\frac{\max S}{\sqrt{2n}} \leq y)=\P^{(d)}_{2n}(\max \e \leq y)+O(n^{1/2})$, 
and the convergence of moments of $\frac{\max S}{\sqrt{2n}}$ under $\P^{(d)}_{2n}$ to 
those of $\max \e$ is also known (see  \cite{FO} and references therein, 
where this is stated in link with the convergence of the height of random trees). 
These convergence results are not strong enough to obtain Theorem \ref{princ}.}
\vskip 0.3cm 
\par 

The strategy will be at first to transform the question in terms of Bernoulli bridges,  
and then to transform the question in terms of simple random walks where the answer is easy. The steps follow some ideas developed by Janson and Marckert \cite{JFM} in their proof of their \mbox{Lemma 1.}

\subsection{From Dyck paths to Bernoulli bridges }

Let us introduce the set $\B_n$  of ``Bernoulli bridges'' with $n$ steps
$$
\B_n=\{ \S: \ \S\in\W_n, S_n=-1\}. 
$$
The quotes around ``Bernoulli bridges'' are there to signal that often the terms ``Bernoulli bridges'' concerns walks ending at 0 instead at $-1$.
 Clearly, $\B_n$ is empty for even $n$ and it is easy to see that $\#\B_{2n+1}=\binom{2n+1}{n}$; we denote by 
 $\P^{(b)}_{2n+1}$ the uniform distribution on 
 $\B_{2n+1}$, and by $\E^{(b)}_{2n+1}$ the expectation with respect to  $\P^{(b)}_{2n+1}$.\par
\vskip 0.2cm 
The cycle Lemma introduced by Dvoretzky and Motzkin \cite{DM}  
(see also Raney \cite{R}  and also Pitman \cite{Pit}, Section 6.1) 
allows one to relate quantities on Dyck paths and on Bernoulli bridges, and among other explains why 
\begin{equation}\label{pfoo}
(2n+1)\#\D_{2n}=\#\B_{2n+1}.
\end{equation}
Consider the set of Dyck paths with size $2n$ with an additional last step $-1$~:
\[\D^{\star}_{2n+1}:=
\{\S~:\S\in\W_{2n+1}, S_i\geq 0 \textrm{ for any }i \in\L{0,2n}\R, S_{2n}=0,S_{2n+1}=-1\}.\]
Obviously there is a canonical correspondence between $\D^{\star}_{2n+1}$ and $\D_{2n}$, and this correspondence conserves the value of the maximum of the paths. Now, the left hand side of (8) is viewed to be the cardinality of 
$\D^{\star}_{2n+1}\times \L{1,2n+1}\R$.\par
 We state the Cycle Lemma as follows:
 
 \vskip 0.3cm
 
\noindent {\bf Lemma 2.}

\noindent  {\it  There exists a 
one-to-one correspondence 
$\Psi_{2n+1}$ between \mbox{$\D^{\star}_{2n+1}\times \L{1,2n+1}\R$} and $\B_{2n+1}$ and such that
if $\S '=\Psi_{2n+1}(\S,k)$ for some $k\in \L{1,2n+1}\R$, then
\begin{equation}\label{borne}
| (\max \S) - (\max \S'-\min \S') | \leq 1.
\end{equation}
}

We provide a proof of this classical result for reader's convenience.\\

\noindent {\bf Proof.} 
For any walk $\S$ in $\W_{n+1}$, let 
\[\Delta_{j}(\S)=\S_{j+1}-\S_{j},~~ j\in\L{0,n}\R \]
denote the list of increments of $\S$.
For a fixed  $(\S,k)$ element in $\D^{\star}_{2n+1}\times \L{1,2n+1}\R$, we let 
$\Psi_{2n+1}(\S,k)$ be the walk whose list of increments is 
$$(\Delta_{(i+k) \mod (2n+1)}(\S), i =0\dots, 2n).$$ 

For any $k$ in  $\L{1,2n+1}\R$,  $\S':=\Psi_{2n+1}(\S,k)$
is indeed a bridge since the sum of the increments is $-1$. We now explain why $\Psi_{2n+1}$ is a bijection from
 $\D^{\star}_{2n+1}\times \L{1,2n+1}\R$ onto $\B_{2n+1}$.  
For a fixed element  $\S \in\D^{\star}_{2n+1} $, let
\[\Psi_{2n+1}(\S):=\left\{\Psi_{2n+1}(\S,k), k \in \L{0,2n}\R \right\},\]  
be ``a rotation class''. It is easy to see that $\Psi_{2n+1}(\S,k)$ reaches its minimum for the first time at time $2n+1-k$. Hence,  $\Psi_{2n+1}(\S,k)\neq\Psi_{2n+1}(\S,k')$ if $k\neq k'$, and each rotation class $\Psi_{2n+1}(\S)$ contains a unique Dyck path. 
\vskip 0.2cm 
It remains to explain why each bridge belongs to a unique rotation class: take a bridge $\S$ that reaches its minimum for the first time at time $k$. The walk $\S'$ whose list of increments is $(\Delta_{(i+2n+1-k) \mod (2n+1)}(\S), i =0\dots, 2n)$ is a Dyck path. Thus,  $\Psi_{2n+1}(\S',k)=\S$, and then $\S$ belongs to the rotation class of $\S'$ (and only to this one).
As a conclusion, each rotation class contains a unique element of $\D^{\star}_{2n+1}$, has cardinality $2n+1$, and of course each element of $\D^{\star}_{2n+1}$ belongs to a rotation class. 
\vskip 0.2cm 
Now it is easy to see that for $\S$ in $\D^{\star}_{2n+1}$ and for any $k\in\L{0,2n}\R$, 
the bridge $\S'=\Psi_{2n+1}(\S,k)$ satisfies inequality $|(\max \S) - (\max \S'-\min \S')| \leq 1$. ~$\Box$\medskip

Hence, the uniform distribution on $\B_{2n+1}$ is the push-forward measure of the uniform distribution on  
$\D^{\star}_{2n+1}\times \L{1,2n+1}\R$ by 
$\Psi_{2n}$ (which amounts to first choosing a Dyck path uniformly, and then a rotation).
It follows from all these considerations that
\[\sup_n \E^{(d)}_{2n}\left(e^{\lambda \frac{\max \S}{\sqrt{2n}}}\right)<+\infty \textrm{ if and only if }
\sup_n \E^{(b)}_{2n+1}\left(e^{\lambda \frac{\max \S-\min \S}{\sqrt{2n}}}\right)<+\infty.\] 
We now show that this second assertion holds.

\subsection{From bridges to simple walks}

For any walk $\S$ let
$$
Y^{\S}_{[a,b]} :=  \max_{a\le k\le b} S_k- \min_{a\le k\le b}S_k.
$$

We have, using 
$
Y^{\S}_{[0,2n+1]} \le Y^{\S}_{[0,n]} + Y^{\S}_{[n,2n+1]}
$ and the Cauchy-Schwarz inequality
\begin{eqnarray}
\E^{(b)}_{2n+1}\left(e^{\lambda \frac{\max \S-\min \S}{\sqrt{2n}}}\right)&=&
\E^{(b)}_{2n+1}\left(e^{\lambda \frac{Y^{\S}_{[0,2n+1]}}{\sqrt{2n}}}\right)\\
&\leq& 
\left[ \E^{(b)}_{2n+1}\left(e^{2\lambda \frac{Y^{\S}_{[0,n]}}{\sqrt{2n}}}\right)
\E^{(b)}_{2n+1}\left(e^{2\lambda\frac{ Y^{\S}_{[n,2n+1]}}{\sqrt{2n}}}\right)  \right]^{1/2}
\end{eqnarray}
The idea here is to work on the half of the trajectory where the conditioning $S_{2n+1}=-1$ will appear to be "not so important".
Since a time reversal of Bernoulli bridges with size $2n+1$, followed by a symmetry with respect to the $x$-axis send $\B_{2n+1}$ onto $\B_{2n+1}$ and exchange the "two halves" of the trajectory, we just have to prove that for $a=0$ and $a=1$ 
\begin{equation}\label{inter}
\sup_n \E^{(b)}_{2n+1}\left(e^{2\lambda Y^{\S}_{[0,n+a]}}\right) < +\infty.
\end{equation}
We will treat the case $a=0$ the other one being very similar. \par
We equip the space $\W_{2n+1}$ with the filtration $\CF:=(\CF_k)$ where $\CF_k$ is generated by the random variables $(S_{1},\dots,S_k)$.

\vskip 0.3cm
\noindent {\bf Lemma 3.} {\it 
  Let $A_n$ be an $\CF_{n}$-measurable event. We have 
  \begin{equation}\label{key}
  \P^{(b)}_{2n+1}(A_n)\leq C_0\P^{(w)}_{2n+1}(A_n)=C_0\P^{(w)}_{n}(A_n),
  \end{equation}
  for a constant $C_0$ valid for all $n$ (and all $A_n$).
}
\vskip 0.3cm

{\bf Proof .} 
The equality in this formula is clear since under $\P^{(w)}_{2n+1}$,
$\S$ is a Markov chain. Only the existence of $C_0$ is needed to be proved. 
In the following computations, we will use that since both $\P^{(w)}_{2n+1}$ and $\P^{(b)}_{2n+1}$ 
are the uniform distributions on their respective set, we have 
$$
\P^{(b)}_{2n+1}=\P^{(w)}_{2n+1}(\ \cdot\  \vert S_{2n+1}=-1).
$$
 We will also use that under both $\P^{(w)}_{2n+1}$ and $\P^{(b)}_{2n+1}$, $\S$ is a Markov chain.\par
If  $A$ is $\CF_n$-measurable, then
$$
\P^{(w)}_{2n+1}(A | S_n=k, S_{2n+1}=-1)=\P^{(w)}_{2n+1}(A|S_n=k)=\P^{(w)}_{n}(A|S_n=k).
$$
This gives   the following chain of  equalities 
$$
\P^{(b)}_{2n+1}(A) = \P^{(w)}_{2n+1}(A | S_{2n+1}=-1)
=\sum_k \P^{(w)}_{n}(A | S_{n}=k) \P^{(w)}_{2n+1}(S_n=k | S_{2n+1}=-1) .
\eqno (14)
$$

Let us denote by 
 $N(m,j)$ 
 the number of trajectories of simple walks of the length $m$ that end at $j$. 
 Clearly
 $N(m,j) = 
{{m}\choose{ {(m+j)/ 2}}},$ when $m$ and $j$ have the same parity, zero if not. 
 Also it is easy to see  that the number of trajectories $\S$  in $\B_{2n+1}$ such that $S_n=j$  
is  $N(n,j)N(n+1,j+1)$. 
Then one obtains by simple counting arguments 
that
\begin{eqnarray*}
\P^{(w)}_{2n+1}(S_n=k | S_{2n+1}=-1)&=&\frac{N(n,k)N(n+1,k+1)}{ N(2n+1,-1)}\\
&=&\frac{N(n,k)}{2^{n}} \frac{2^n N(n+1,k+1)}{N(2n+1,-1)}\\
&=& \P^{(w)}_{n}(S_n=k) \frac{2^n N(n+1,k+1)}{N(2n+1,-1)}
\leq  C_0 \P^{(w)}_{n}(S_n=k),
\end{eqnarray*}
where
$$
C_0 = \sup_{n\ge 1} \sup_{k}\frac{2^n N(n,k+1)}{N(2n+1,-1)} = 
\sup_{n\ge 1}  2^n { n\choose \L n/2\R} {2n\choose n}^{-1},
$$
is indeed finite (as one may check using the Stirling formula). 
Hence, the right hand side in (14) is bounded by 
\[\sum_k \P^{(w)}_{n}(A | S_{n}=k)\cdot C_0 \P^{(w)}_{n}(S_n=k)=C_0\P^{(w)}_{n}(A ).\]
This ends the proof of the Lemma. $\Box$
\vskip 0.3cm 

To conclude the proof of Theorem 1, we explain why (12) holds true.  Using Lemma 3, we have
$$
  \sup_n \E^{(b)}_{2n+1}\left(e^{2\lambda Y^{\S}_{[0,n]}}\right) \leq C_0 \sup_n 
  \E^{(w)}_{n}\left(e^{2\lambda Y^{\S}_{[0,n]}}\right) ,
 \eqno (15)
 $$ 
  Since $Y^{\S}_{[0,n]}$ is $\CF_{n}$-measurable.
The right hand side of (15) is much simpler than the left one, since it deals with simple random walks under the uniform distribution. Then using again Cauchy-Schwarz, it suffices to show that
$$
 \sup_n \E^{(w)}_{n}\left(e^{4\lambda\frac{\max_{0\leq k \leq n} S_k}{\sqrt{2n}}}\right)<+\infty
\eqno (16)
$$
and the same thing for $\max$ replaced by $\min$ (which gives the same quantity).
Now, by the Andr\'e reflexion principle (see Feller \cite{FEL}, page 72), we have 
\[\P^{(w)}_{n}(\max_{0\leq k\leq n} S_k \geq x)=2\P^{(w)}_{n}(S_n>x)+\P^{(w)}{n}(S_n=x)\leq 2\P^{(w)}_{n}
(S_n\geq x).\]
Now, the use of Hoeffding's inequality yields directly (16).


\begin{thebibliography}{99}

\bibitem{BIL} P. Billingsley. {\em Probability and measure.} Second Edition, {\it John Wiley \& Sons, Inc., New York } (1986).
MR0830424

\bibitem{CHU} K.L. Chung. {\em Maxima in Brownian excursions.}  Bull. Amer. Math. Soc.  81  (1975), 742--745. MR0373035

\bibitem{DM}  A. Dvoretzky and Th. Motzkin. {A problem of arrangements.} 
Duke Math. J. {\bf 14}\rm  (1947), 305- 313. 
MR0021531

\bibitem{FEL} W. Feller.  {\em An introduction to probability theory and its applications.} Vol. I. Third edition, John Wiley \& Sons, Inc., New York-London-Sydney (1968). 
MR0228020


\bibitem{FO} Ph. Flajolet, A. Odlyzko.  The average height of binary trees and other simple trees. 
{\it J. Comput. System Sci.,} {\bf 25} (1982),  171--213. 
MR0680517

\bibitem{JFM} S. Janson, J.-F. Marckert.  Convergence of Discrete Snakes.  {\it J.  Theor. Prob.}
 {\bf 18}   (2005), 615-- 645. MR2167644


\bibitem{KAI} W.D. Kaigh. {\ An invariance principle for random walk conditioned by a late return to zero.} 
{\it Ann. Probab.} {\bf 4} (1976), 115-121. MR0415706

\bibitem{KV} O. Khorunzhiy, V. Vengerovsky. { Ewen walks and estimates of high moments of large Wigner random matrices}. {\it Preprint arXiv:0806.0157.}  Math. Review number not available



\bibitem{R} G. N. Raney. { Functional composition patterns and power series reversion.}
{\it Trans. Amer. Math. Soc.} {\bf 94} (1960), 441-451. MR0114765

\bibitem{STAN} R.P. Stanley. {\em Enumerative combinatorics}. Vol. 2, Cambridge University Press. (1999) 
MR1676282



\bibitem{Pit} J. Pitman. {\em  Combinatorial stochastic processes.} Lecture Notes in Mathematics Vol. 1875, Springer-Verlag, Berlin (2006). Lectures from the 32nd Summer School on Probability Theory held in Saint-Flour, July 7--24, 2002. MR2245368

\bibitem{SS} Ya. G. Sinai, A. Soshnikov. { A refinement of Wigner's semicircle law in a neighborhood of the spectrum edge for random symmetric matrices}. {\it Func. Anal. Appl.} {\bf 32} (1998), 114-131. MR1647832

\bibitem{SD} L. Smith,   P. Diaconis. {   Honest {B}ernoulli excursions.}
{J. Appl. Probab.}  {\bf 25}  (1988), 464--477. MR0954495

\bibitem{S} A. Soshnikov. {Universality of the edge of the spectrum in Wigner random matrices.} 
{\it Commun. Math. Phys.} {\bf 207} (1999), 697-733. 
MR1727234

\bibitem{W} E. Wigner. { Characteristic vectors of bordered matrices with infinite dimensions}.
 {\it Ann. Math.} {\bf 62} (1955), 548-564. 
MR0077805


\end{thebibliography}
\end{document}